

%




\magnification=\magstep1
\baselineskip 18pt
\def\:{\Vert}
\def\({\lbrack}
\def\){\rbrack}
\def\S{{\cal S}}
\def\N{{\cal N}}
\def\M{{\cal M}}
\def\T{{\bf T}}
\def\R{{\bf R}}
\def\A{{\cal A}}
\def\B{{\cal B}}
\def\S{{\cal S}}

\def\endpf{\ \ \vrule height6pt width4pt depth2pt}


\centerline{\bf NEST ALGEBRAS IN $c_1$}

\centerline{by}

\centerline
{Alvaro Arias\footnote{*}{Research partially supported by  
BSF89--00087}}

\centerline{Department of Theoretical Mathematics}

\centerline{The Weizmann Institute of Science}

\centerline{Rehovot, Israel}
\bigbreak

ABSTRACT. In this paper we address some basic questions of the
Banach space structure of the nest algebras in the trace class; in
particular, we study whether any two of them are isomorphic to each other,
and show that the nest algebras in the trace class have bases. We
construct three non-isomorphic examples of nest algebras in $c_1$;
present a new proof of the
primarity of $c_1$ (Arazy, [Ar1], [Ar2]), and prove that $K(H)$, and the
nest algebras in $B(H)$ are primary.
\bigbreak

\noindent{{1. INTRODUCTION.}}
\medbreak
In the present paper we study some basic questions of
the Banach space structure of the nest algebras.
In particular, we study whether any two nest algebras in $c_1$
are isomorphic to each other.

The answer to this question is known for the other Schatten
$p$-classes, $c_p$, and for $B(H)$: All the nest algebras in
$c_p$, $1<p<\infty$ are isomorphic to $c_p$. This is an easy
consequence of the results of Macaev \(Ma\) and
Gohberg and Krein \(GK\) that say that the nest algebras
in $c_p$, $1<p<\infty$ are complemented in $c_p$. Likewise,
all the nest algebras in $B(H)$ are completely isomorphic to
each other (see  \(A2\)).

The structure of the nest algebras in $c_1$ is richer.
We will show, for instance, that
if the complete nest is uncountable, then the nest algebra
in $c_1$ is isomorphic to the continuous nest; and for the countable  
case
there is a natural collection of spaces, indexed by the countable
ordinal numbers, that resembles the classification of spaces of
continuous functions on countable metric spaces given by Bessaga and
Pe\l czy\'nski \(BP\); although we can only prove that three of them
are not isomorphic to each other.

The Banach space invariant we use is primarity, (a Banach space $X$
is primary if whenever $X\approx Y\oplus Z$ then either $X\approx Y$
or $X\approx Z$). As  side results we prove that the nest algebras in
the trace class have bases and that the trace class,
the space of compact operators and the nest algebras in $B(H)$ are
primary. J. Arazy, [Ar1], [Ar2], gave an earlier proof (1980-1) of
the primarity of $c_1$.
Our proof is shorter and extends to $K(H)$;
the technique was motivated by a paper of Blower
[Bl].

\medbreak
Section 2 has the preliminaries; we fix the notation
and quote the necessary results from operator and Banach space
theory needed later.  In Section 3 we prove that $c_1$ is
primary. In Section 4 we
apply the technique developed in Section 3 to study the 
nest algebras in $c_1$.  In Section 5 we prove that
the nest algebras in $B(H)$ are primary. We conclude with
applications and open questions in Section 6; in particular, we
prove that the nest algebras in the trace class have bases.
\medbreak
The author wants to thank G. Schechtman for the great hospitality during
his year at the Weizmann Institute of Science and to
J. Arazy for explaining the content of \(Ar2\).

\bigbreak
\noindent{{2. PRELIMINARIES.}}
\medbreak
For this paper $H$ denotes a separable Hilbert space and
$B(H)$ the set of all linear, bounded operators on $H$. 
A complete {\it nest}, $\N$, is a totally ordered family of closed
subspaces that contains $0$, $H$, and is closed under intersections
and closed unions. The {\it nest algebra} induced by $\N$ is the
set of all $T\in B(H)$ that leave invariant the elements of $\N$;  
i.e.,
$$\hbox{Alg}\N=\{T\in B(H):TN\subseteq N\hbox{ for every }N\in\N\}.$$

The following examples have motivated a big part of the theory.
\medbreak
{EXAMPLE 1.}
In $\ell_2$ let $N_k=\overline{\hbox{span}}\{e_i\}_1^k$ and
$k=1,2,\cdots,\infty$. Then $\N=\{0\}\bigcup \{N_k\}_1^\infty$ is
a nest and it is easy to see that Alg$\N$ is the
set of upper triangular operators.
\medbreak
{EXAMPLE 2.} In $L_2(0,1)$ let $N_t=\{f\in L_2:
\hbox{supp} f\subset\(0,t\)\}$
for $0\leq t\leq1$.
\medbreak
{EXAMPLE 3.} In $L_2(0,1)\oplus L_2(0,1)$ let
$M_t=\{f\in L_2\oplus L_2:
\hbox{supp} f\subset \(0,t\)\oplus\(0,t\)\}$
for $0\leq t\leq1$.
\medbreak
The nest algebras were introduced by Kadison and Singer'60 \(KS\)
and Ringrose'65 \(R\). A central problem was that of classification.
Two nests $\N$ and $\M$ are {\it similar} (unitarily equivalent)
if there exists $T\in B(H)$ {\it invertible} (unitary) such that
$T\N=\M$; equivalently, $T\hbox{Alg}\N T^{-1}=\hbox{Alg}\M$.

It is clear, (by a simple cardinality argument), that the nests
of examples 1 and 2 are not similar; however, it was open for
a long time whether the nests of examples 2 and 3 were. This was  
answered
by Larson'85 \(L\), who not only proved they are similar but also
that any two continuous nests (i.e., those whose index
is connected for the order topology) are similar.

Examples 2 and 3 are particular cases of the following
natural family of nest algebras:
Let $L_2(\mu)=L_2(\(0,1\),\mu)$ where $\mu$ is a positive Borel
measure on $\(0,1\)$. For $0\leq t\leq1$, let
$M_t=\{f\in L_2(\mu):\hbox{supp}f\subseteq\(0,t\)\}$, and
$M_t^-=\{f\in L_2(\mu):\hbox{supp}f\subseteq\(0,t)\}$.  Then
$\M=\{M_t,M_t^-\}_{0\leq t\leq1}$ is a complete nest called the  
standard
nest.

J. Erdos \(E\)
proved that if $\N$ is a nest in a separable Hilbert space,
then there exists a sequence $\mu_1>>\mu_2>>\cdots$ of regular Borel
measures on $\(0,1\)$ such that Alg$\N$ is unitarily equivalent to
the standard nest on
$L_2(\mu_1)\oplus L_2(\mu_2)\oplus\cdots$; i.e., for $0\leq t\leq1$
and $f\in L_2(\mu_1)\oplus L_2(\mu_2)\oplus\cdots$,
we have that
$f\in M_t$ if and only if  
supp$f\subset\(0,t\)\oplus\(0,t\)\oplus\cdots$,
and
$f\in M_t^-$ if and only if
supp$f\subset\(0,t)\oplus\(0,t)\oplus\cdots$.

\medbreak
For $1\leq p<\infty$ let $c_p$, the Schatten $p$-class, be the set of  
all
$T\in B(H)$ for which $\|T\|_p^p=\hbox{tr}(T^*T)^{p/2}$ is finite. Given
$\N$, a nest in $H$, define 
$(\hbox{Alg}\N)^p=\hbox{Alg}\N\bigcap c_p$
to be the corresponding nest algebra in $c_p$.
Macaev \(M\) and Gohberg and Krein \(GK\)
proved that an infinite nest algebra in $c_p$ is complemented in
$c_p$ if and only if $1<p<\infty$. Since this behavior 
is identical to that of the Hardy spaces $H^p$
in $L_p$, the nest algebras in $c_p$ are sometimes called
the non-commutative $H^p$-spaces. However, there are many more
analogies than that (see for example \(FAM\), \(P\), \(A1\)).

The similarity theorem also extends to the nest algebras in $c_p$.
If $\N$
and $\M$ are continuous nests then we can find $T\in B(H)$ invertible
such that $T\hbox{Alg}\N T^{-1}=\hbox{Alg}\M$.
Which implies, of course, that Alg$\N\approx\hbox{Alg}\M$.
But since $T,T^{-1}\in B(H)$ we have that
$T(\hbox{Alg}\N)^p T^{-1}=(\hbox{Alg}\M)^p$; hence,
$(\hbox{Alg}\N)^p\approx(\hbox{Alg}\M)^p$.
In particular, up to similarity, there is only one continuous
nest in $c_1$ which we
denote by $\T^1(\R)$.
\medbreak

Example 1 has been studied from the Banach space point of view
where it is denoted by $\T$ (triangular) and $\T^p=\T\bigcap c_p$.
It is a particular case of \(GK\) that $\T^p$ is complemented in  
$c_p$
if and only if $1<p<\infty$. This fact was stressed by Arazy \(Ar2\)  
who
proved that $\T^1$ is not isomorphic to a complemented subspace of  
$c_1$, we will use this fact in the proof of Proposition 11.
Another important fact, proved by Kwapien and Pe\l czy\'nski \(KP\),
says that $c_1$ does not embed into $\T^1$.
\medbreak

We will use repeatedly the Pe\l czy\'nski decomposition method \(Pe\).  
The form we use asserts that if $X$ embeds complementably into $Y$, $Y$ embeds
complementably into $X$ and $X\approx(\sum\oplus X)_p$ for some
$1\leq p\leq\infty$ then $X\approx Y$.

A Banach space $X$ is primary if whenever $X\approx Y\oplus Z$
then either $X\approx Y$ or $X\approx Z$.
It is an immediate consequence of Pe\l czy\'nski's decomposition
method that if $X\approx(\sum\oplus X)_p$ for $1\leq p\leq\infty$  
then
$X$ is primary if for any $T:X\to X$ bounded and linear,
the identity on $X$ factors through $T$ or through $I-T$.
($I$ factors through $T$ if we can find $A,B:X\to X$ bounded and
linear for which $I=ATB$).

Finally, the necessary combinatorial
results used in Section 5 can be found in \(Bo\).
\medbreak

\bigbreak
\noindent{3. $c_1$ IS PRIMARY}
\medbreak
In this section we give a new proof of the primarity of $c_1$. The technique
of the proof will be used in the next section to distinguish
different isomorphic types of nest algebras in $c_1$.
\medbreak
Let $(e_i)_{i=1}^\infty$ be an orthonormal basis for $H$, and
let $e_{ij}=e_j\otimes e_i$
be the rank-1 operator  sending $z$ to $(z,e_j)e_i$.
Let $\sigma$ and $\psi$ be infinite subsets of {\bf N} and 
define $J_{\sigma,\psi}:c_1\to c_1$ and
$K_{\sigma,\psi}:c_1\to c_1$ by
$$\eqalign{J_{\sigma,\psi}(e_{ij})&=e_{\sigma(i)\psi(j)},\hbox{  
and}\cr
           K_{\sigma,\psi}(e_{ij})&=
     \cases{e_{kl},&if $\sigma(k)=i\hbox{ and }\psi(l)=j$;\cr
              0   ,&otherwise.\cr}\cr}$$
These maps were used in \(KP\); however, our notation
and motivation comes from a paper of Blower \(Bl\) where he proved a finite
dimensional analogue of Theorem 1 below.

We will use $\sigma$ and $\psi$
as subsets or as functions $\sigma:{\bf N}\to{\bf N}$ according
to our needs; $\sigma(i)$ denotes the $i$th smallest element
of $\sigma$.
It is easy to see that $J_{\sigma,\psi}$
is an isometric embedding and $K_{\sigma,\psi}J_{\sigma,\psi}=I$,
$I$ is the identity of $c_1.$
Moreover, if $\sigma_i$ and $\psi_i$, $i=1,2$ are infinite
subsets of ${\bf N}$ then $J_{\sigma_1,\psi_1}J_{\sigma_2,\psi_2}=
J_{\sigma_1\sigma_2,\psi_1\psi_2}$, where
$\sigma_1\sigma_2(j)=\sigma_1(\sigma_2(j))$. 
Similarly, $K_{\sigma_1,\psi_1}K_{\sigma_2,\psi_2}=
K_{\sigma_1\sigma_2,\psi_1\psi_2}$.

For this section, $\Phi$, with or without subscripts, denotes a bounded linear
operator.

\medbreak
THEOREM 1. For every $\epsilon>0$ and
$\Phi:c_1\to c_1$, we can find
$\sigma,\psi\subset{\bf N}$, and
$\lambda\in{\bf C}$ such that
$\:K_{\sigma,\psi}\Phi J_{\sigma,\psi}-\lambda I\:<\epsilon.$
Thus, one of $K_{\sigma,\psi}\Phi J_{\sigma,\psi}$ and
$I-K_{\sigma,\psi}\Phi J_{\sigma,\psi}$ is invertible.
\medbreak

We prove Theorem 1 in 5 steps. Each one of them will be
a factorization of the form $K_{\sigma,\psi}\Phi J_{\sigma,\psi}$;
yielding new
maps $\Phi_i$ with nicer properties. The sets $\sigma$ and $\psi$ are
constructed inductively.

REMARK. It might be instructive to consider the following example ``far''
from a multiplier. Let $\Phi:c_1\to c_1$ be the transpose operator; i.e.,
$\Phi e_{ij}=e_{ji}$; $\sigma$ the set of even integers and $\psi$ the set of
odd integers. Then $K_{\sigma,\psi}\Phi J_{\sigma,\psi}=0$.
\smallbreak

We will use several times the
following elementary lemma.
\medbreak
{LEMMA 2.} Let $N$ be either finite or infinite,
$E_n$ an $n$-dimensional space,
$\epsilon>0$ and $T:\ell_2^N\to E_n$ a bounded, linear map. Then,
card$\{i\leq N:\|Te_i\|>\epsilon\}\leq n^3\|T\|^2/\epsilon^2.$
\medbreak
PROOF. Let $\{\tilde{e}_i\}_{i\leq n}$ be an Auberbach basis
for $E_n$; i.e.,
$$\max_{i\leq n}\vert a_i\vert\leq
\bigl\|\sum_{i\leq n}a_i\tilde{e}_i\bigr\|
\leq\sum_{i\leq n}\vert a_i\vert.$$
For every $i\leq n$ let $A_i=\{j\leq N: |\tilde{e_i}^*(T\tilde{e}_j)|
>\epsilon/n\}$, where $\{\tilde{e_i}^*\}_{i\leq n}$ is the dual basis in $E_n^*$.
Choose $\vert\epsilon_j\vert=1$ appropriately so that
$\|\sum_{j\in A_i}\epsilon_je_j\|=\sqrt{\hbox{card}A_i}$ and
$\|T(\sum_{j\in A_i}\epsilon_je_j)\|\geq\epsilon\hbox{ card}A_i/n$.
Then it is clear that card$(A_i)\leq n^2\|T\|^2/\epsilon^2.$
Hence, if $j\not\in\bigcup_{i\leq n}A_i$ we have that
$\|Te_j\|<\epsilon$ and card$(\bigcup_{i\leq n}A_i)\leq n^3\|T\|^2
/\epsilon^2$\endpf
\medbreak
REMARK. We will not need the estimate of Lemma 2. It suffices to know that
card$\{i\leq N:\|Te_i\|>\epsilon\}$
is small compared to $N$ and independent of $N$. Moreover, we will also
apply Lemma 2
for $T:X\to E_n$ when $X$ is just isomorphic, not
necessarily isometric, to a Hilbert space. The  
estimate changes but depends on the Banach-Mazur distance of $X$ to 
the respective Hilbert space, and not on the dimension of $X$.

\medbreak

PROOF OF THEOREM 1.
We introduce some notation now. For every $n\in{\bf N}$ let
$F_n={\rm span }\{e_{ij}:\max\{i,j\}=n\}$ and
$H_n=\overline{\rm span }\{e_{ij}:\min\{i,j\}=n\}.$ Notice that
both $\{F_n\}$ and $\{H_n\}$ form a Schauder decomposition for
$c_1$. The first one has very nice properties (see \(KP\) and
\(AL\) ), and $H_n\approx\ell_2$.

We also use $M_n={\rm span }\{e_{ij}:\max\{i,j\}\leq n\}$ with
$P_n$ the natural projection onto it; and
$E_n=\overline{\rm span }\{e_{ij}:\min\{i,j\}\leq n\},$ with
$Q_n$ the natural projection onto it. $E_n$ is still isomorphic
to a Hilbert space (but with an isomorphism constant depending on  
$n$).

\medbreak
STEP 1. For every $\epsilon>0$, there exist $\sigma_1\subset{\bf N}$ 
and $\Phi_1$ such that
$\|\Phi_1-K_{\sigma_1,\sigma_1}\Phi  
J_{\sigma_1,\sigma_1}\|<\epsilon$, and
$\Phi_1M_n\subset M_n$, $\Phi_1E_n\subset E_n$ for every $n\in{\bf N}$.
\medbreak
The proof of Step 1 is easy. We present in full detail the  
construction for
$M_n$ and indicate how to do it for $E_n$. 

The key ideas are that if $K\subset c_1$ is compact, then
there is some $m$ such that $K$ is essentially inside $M_n$; and 
if $E\approx\ell_2$, then there is some $m$ such that $E$ is  
essentially
inside $E_m$. (See [Ar1], Proposition 2.2).

Let $\sigma_1(1)=1$ and assume that we have chosen
$\sigma_1(1),\cdots,\sigma_1(n)$. 
Since $\Phi\, Ball(M_{\sigma_1(n)})$
is compact we can find $m>\sigma_1(n)$ such that 
$\sup_{x\in\, Ball(M_{\sigma_1(n)})}\|P_m\Phi x-\Phi  
x\|<\epsilon_{n+1}$,
where $\epsilon_{n+1}>0$ is chosen small enough. Then set 
$\sigma_1(n+1)=m+1.$ Proceeding this way we construct $\sigma_1$.

Let $x\in F_n$. Then $J_{\sigma_1,\sigma_1}x\in F_{\sigma_1(n)}$; and
hence, $\|\Phi J_{\sigma_1,\sigma_1}x-P_m\Phi  
J_{\sigma_1,\sigma_1}x\|\leq
\epsilon_n\|x\|$, where $m=\sigma_1(n+1)-1$.  It is easy to check  
that
$K_{\sigma_1,\sigma_1}P_m=P_{\sigma_1(n)}K_{\sigma_1,\sigma_1}$.  
Therefore,
$$\|K_{\sigma_1,\sigma_1}\Phi J_{\sigma_1,\sigma_1}x-P_{\sigma_1(n)}
K_{\sigma_1,\sigma_1}\Phi J_{\sigma_1,\sigma_1}x\|<\epsilon_n\|x\|.$$

Define $\Phi'$ by $\Phi'x=P_{\sigma_1(n)}K_{\sigma_1,\sigma_1}\Phi
J_{\sigma_1,\sigma_1}x$ for $x\in F_n$. Then, if  
$\sum_n\epsilon_n<\epsilon$
is small enough, $\Phi'$ is well defined and satisfies
$\|K_{\sigma_1,\sigma_1}\Phi J_{\sigma_1,\sigma_1}-\Phi'\|<\epsilon$  
and
$\Phi' M_n\subset M_n$ for every $n\in{\bf N}$.

Repeat the process for $\Phi'$ with respect to the $E_n$'s, doing the
perturbation argument along the $H_n$'s and finish.

\medbreak
STEP 2. For every $\epsilon>0$ there exit $\sigma_2,\psi_2\subset{\bf  
N}$ and $\Phi_2\in B(c_1)$ such that
$\|\Phi_2-K_{\sigma_2,\psi_2}\Phi_1 J_{\sigma_2,\psi_2}\|<\epsilon$,  
and
$\Phi_2F_n\subset F_n$, $\Phi_2 E_n\subset E_n$ for every $n\in{\bf  
N}$.
\medbreak

We construct $\sigma_2,\psi_2$ satisfying
$$\sigma_2(1)\leq\psi_2(1)<\sigma_2(2)\leq\psi_2(2)<\sigma_2(3)\leq
\psi_2(3)
\cdots,\leqno{(1)}$$
to guarantee that $M_n$ and $E_n$ are invariant for $\Phi_2$.

Let $\sigma_2(1)=\psi_2(1)=1$ and assume that we have chosen
$\{\sigma_2(1),\cdots,\sigma_2(n)\}$ and
$\{\psi_2(1),\cdots,\psi_2(n)\}$ satisfying
$\sigma_2(1)\leq\psi_2(1)<\cdots<\sigma_2(n)\leq\psi_2(n)$.

Let $N>\psi_2(n)$ be a ``large'' number and consider
$$A=\{j>N\colon\hbox{ for }1\leq i\leq N,\,\|P_{\psi_2(n)}\Phi_1
e_{i,j}\|<\epsilon_{n+1}\}.$$ 

By Lemma 2, $A^c$ is finite. Choose
$\psi_2(n+1)=\min A$, and since $N$ is large enough, can  
find $i_0$, $\psi_2(n)<i_0\leq N$ satisfying: for all $1\leq j\leq n$,
$\|P_{\psi_2(n)}\Phi_1 e_{i_0\psi_2(j)}\|\leq\epsilon_{n+1}$. Then set 
$\sigma_2(n+1)=i_0$. Proceeding this way we construct  
$\sigma_2,\psi_2$.

Summarizing we have: If $e_{ij}\in F_n$, then
$\|P_{\psi_2(n-1)}\Phi_1 J_{\sigma_2,\psi_2}e_{ij}\|<\epsilon_n$. It  
is
easy to see that (1) implies  
$K_{\sigma_2,\psi_2}P_{\psi_2(n-1)}=P_{n-1}K_{\sigma_2,\psi_2}$;  
hence,
$\|P_{n-1}K_{\sigma_2,\psi_2}\Phi_1  
J_{\sigma_2,\psi_2}e_{ij}\|<\epsilon_n$.

Define  
$\Phi_2e_{ij}=(P_n-P_{n-1})K_{\sigma_2,\psi_2}\Phi_1J_{\sigma_2,\psi_ 
2}
e_{ij}$, where $e_{ij}\in F_n$ (i.e., the projection onto $F_n$.  
Recall that $K_{\sigma_2,\psi_2}\Phi_1J_{\sigma_2,\psi_2}M_n\subset  
M_n$).
Therefore, 
$$\|(K_{\sigma_2,\psi_2}\Phi_1J_{\sigma_2,\psi_2}-\Phi_2)e_{ij}\|<
\epsilon_n.$$
Since the $(2n-1)$-dimensional space $F_n$ has a 1-basis consisting  
of $e_{ij}$'s, we conclude that if $x\in F_n$, then
$\|(K_{\sigma_2,\psi_2}\Phi_1J_{\sigma_2,\psi_2}-\Phi_2)x\|<
(2n-1)\epsilon_n\|x\|.$

If we choose $\sum_n(2n-1)\epsilon_n<\epsilon$ small enough we  
finish.

\medbreak
STEP 3. For every $\epsilon>0$ there exist $\sigma_3,  
\psi_3\subset{\bf N}$,
and $\Phi_3\in B(c_1)$  such that 
$\|\Phi_3-K_{\sigma_3,\psi_3}\Phi_2 J_{\sigma_3,\psi_3}\|<\epsilon$,
and $\Phi_3H_n\subset H_n$ and $\Phi_3 F_n\subset F_n$
for every $n\in{\bf N}$.
\medbreak

We will choose $\sigma_3,\psi_3$ as in (1) to guarantee that $F_n$ is  
invariant for $\Phi_3$.
Let $\sigma_3(1)=\psi_3(1)=1$ and $A(1)=B(1)={\bf N}$.
Assume that we have chosen $\sigma_3(1),
\cdots,\sigma_3(n)$; $\psi_3(1),\cdots,\psi_3(n)$ and $A(n),B(n)$,
infinite subsets of ${\bf N}$, satisfying:
$\sigma_3(1)\leq\psi_3(1)<\cdots<\sigma_3(n)\leq\psi_3(n)$.
(We will choose $\sigma_3(n+1)$ from
$A(n)$ and $\psi_3(n+1)$ from $B(n)$).

Let $N>\psi_3(n)$ be a ``large'' number. For every $j\in B(n)$,
$j>N$ find $i(j)$, $\psi_3(n)<i(j)<N$ such that 
$\|Q_{\psi_3(n)}\Phi_2e_{i(j)j}\|<\epsilon_{n+1}$,
(Apply Lemma 2 to $Q_{\psi_3(n)}\Phi_2:F_j\to Q_{\psi_3(n)}F_j$,
notice that dim$\,(Q_{\psi_3(n)}F_j)=2\psi_3(n)$).
Let $B(n+1)\subset B(n)$
be an infinite subset of those $j$'s with common $i(j)=i_0$ and set
$\sigma_3(n+1)=i_0$. Exchanging the roles of $A(n), B(n)$ with
$B(n), A(n+1)$ we find $A(n+1)\subset A(n)$, and $\psi_3(n+1)$ such  
that
$\|Q_{\psi_3(n)}\Phi_2 e_{i\psi_3(n+1)}\|<\epsilon_{n+1}$ 
for every $i\in A(n+1)$. Proceeding this way we construct  
$\sigma_3,\psi_3$.

Summarizing we have: If $e_{ij}\in H_n$, then  
$\|Q_{\psi_3(n-1)}\Phi_2
J_{\sigma_3,\psi_3}e_{ij}\|<\epsilon_n$. It is easy to check that
$K_{\sigma_3,\psi_3}P_{\psi_3(n-1)}=P_{n-1}K_{\sigma_3,\psi_3}$;  
hence,
$\|Q_{n-1}K_{\sigma_3,\psi_3}\Phi_2J_{\sigma_3,\psi_3}e_{ij}\|<\epsilon_n$.

Define $\Phi_3e_{ij}=(Q_n-Q_{n-1})K_{\sigma_3,\psi_3}\Phi_2
J_{\sigma_3,\psi_3}e_{ij}$, where $e_{ij}\in H_n$ (i.e., the  
projection
onto $H_n$. Recall that  
$K_{\sigma_3,\psi_3}\Phi_2J_{\sigma_3,\psi_3}
E_n\subset E_n$). Therefore, 
$\|(K_{\sigma_3,\psi_3}\Phi_2
J_{\sigma_3,\psi_3}-\Phi_3)e_{ij}\|<\epsilon_n$.
Since $E_n$ is $K_n$-isomorphic to $\ell_2$, and
$Q_nK_{\sigma_3,\psi_3}\Phi_3J_{\sigma_3,\psi_3}:H_{n+1}\to E_n$ 
is ``diagonal'' with respect to the
decompositions: $(H_n\bigcap F_j)_j$ for $H_n$, and $(E_n\bigcap  
F_j)_j$
for $E_n$; we see that if $x\in H_n$, then
$\|(K_{\sigma_3,\psi_3}\Phi_2J_{\sigma_3,\psi_3}-
\Phi_3)x\|<K_n\epsilon_n\|x\|$. 

Since the $H_n$'s form
a Schauder decomposition, it is enough to choose 
$\sum_nK_n\epsilon_n<\epsilon$ small enough to finish.

\medbreak
STEP 4. Find $\sigma_4, \psi_4$ such that $\Phi_4=
K_{\sigma_4,\psi_4}\Phi_3 J_{\sigma_4,\psi_4}$ satisfies
$\Phi_4e_{ij}=\lambda_{ij}e_{ij}$ for some $\lambda_{ij}\in
{\bf C}.$
\medbreak

Just take $\sigma_4=\{1,3,5,\cdots\}$ and
$\psi_4=\{2,4,6,\cdots\}.$  To see that it suffices, notice that
if $i<n$ then $H_i\bigcap F_n=\(e_{in},e_{ni}\);$ hence,
$\Phi_3e_{in}=c_1e_{in}+c_2e_{ni}$ for some constants $c_1,c_2.$

\medbreak
STEP 5. For every $\epsilon>0$, there exist  
$\sigma_5,\psi_5\subset{\bf N}$;
$\Phi_5\in B(c_1)$; and $\lambda\in{\bf C}$ such that
$\|\Phi_5-K_{\sigma_5,\psi_5}\Phi_4 J_{\sigma_5,\psi_5}\|<\epsilon$,
and $\Phi_5e_{ij}=\lambda e_{ij}$ for every $i,j\in{\bf N}$.
\medbreak

Look at the upper part of $\{\lambda_{ij}\}_{i<j}.$
By a standard diagonal argument we find a subsequence
$\sigma_5(1)<\psi_5(1)<\sigma_5(2)<\psi_5(2)<\sigma_5(3)<\psi_5(3)
<\cdots$
such that for some $\lambda_i$, $\lambda\in{\bf C}$,
$\vert \lambda_{\sigma_5(i),\psi_5(j)}-\lambda_i\vert<
\epsilon/2^{i+j}$ and $\vert \lambda_i-\lambda\vert<
\epsilon/2^i$.  Which roughly speaking says that the upper
triangular part of $\tilde{\Phi}_5$ is essentially $\lambda$.
We order them as in (1) to preserve the upper triangular structure
of $\{\lambda_{ij}\}$.
Do the same for the lower part, and assume that it is ``essentially''
$\mu\in{\bf C}$.

Define $\tilde{\Phi}_5=K_{\sigma_5,\psi_5}\Phi_4J_{\sigma_5,\psi_5}$.
Hence, $\tilde{\Phi}_5$
has essentially upper triangular
part $\lambda$ and lower triangular part $\mu$.
Since $\tilde{\Phi}_5$ is bounded they {\it must} agree. Otherwise,
$(\tilde{\Phi}_5-\mu I)/(\lambda-\mu)$ would be like the upper
triangular projection; and the latter one is known to be
unbounded (see \(GK\)).

Let $\Phi_5=\lambda I$, and notice that if $e_{ij}\in H_n$, then
$\|\tilde{\Phi}_5e_{ij}-\Phi_5e_{ij}\|<\epsilon({1\over  
2^{i+j}}+{1\over 2^i})$.
Moreover, it is clear that if $x\in H_n$, then
$\|\Phi_5 x-\tilde{\Phi}_5 x\|<{\epsilon\over 2^{n-1}}\|x\|$. Since
the $H_n$'s form a Schauder decomposition, we finish.
\endpf
\medbreak
COROLLARY 3. (J. Arazy)  $c_1$ is primary.
\medbreak
PROOF. It follows from Theorem 1 that $I$, the identity on  
$c_1$,
factors through $\Phi$ or through $I-\Phi$. This implies 
that if $c_1\approx X\oplus Y$ then $c_1$ embeds complementably
into $X$ or $Y$. Since $c_1\approx(\sum\oplus c_1)_1$, the Pe\l czy\'nski
decomposition method gives the result. \endpf
\medbreak
Notice that $J_{\sigma,\psi}$ and $K_{\sigma,\psi}$ can be defined
in $K(H)$, the space of compact operators in the Hilbert space $H$.
Moreover, it is easy to see that $J_{\sigma,\psi}^*=K_{\sigma,\psi}$,
and $K_{\sigma,\psi}^*=J_{\sigma,\psi}$. Hence, one gets
the equivalence of Theorem 1 and,
\medbreak
COROLLARY 4. $K(H)$ is primary.
\medbreak
REMARKS. (1)
The proof of Theorem 1 works in more general situations. For
instance, if $\Phi:{\bf T^1}\to
{\bf T^1}$ and we make sure that all of the
$\sigma_i$, $\psi_i$, $i=1,\cdots,5$ respect triangularity,
(i.e., they satisfy (1)), then
the same result holds; giving another proof of the fact that
${\bf T^1}$ is primary. It also works for ${\bf T^p}$,
$1<p<2$ giving the same conclusion (both results are proved by
J. Arazy \(Ar1\) ). And for ${\bf T}_E$ if $E$ is a 1-symmetric
sequence space of type $p$, $p<2$.
\medbreak
(2) Some steps of the proof can be adapted to more general
subspaces $\S \subset c_1$ provided we can find enough
$\sigma, \psi\subset{\bf N}$ satisfying $J_{\sigma,\psi}\S
\subset\S$ and $K_{\sigma,\psi}\S\subset\S.$ This fact will be
essential in the next section.

\bigbreak
\noindent{4. NEST ALGEBRAS IN $c_1$}
\medbreak
In this section we study the isomorphism types of the
nest algebras in $c_1$.
Notice that $(\N,\leq)$ is a compact space with the order topology.

The results we obtain are:
\medbreak
{THEOREM 5.} If $\N$ is an uncountable nest then
$(\hbox{Alg}\N)^1$ is isomorphic to $\T^1(\R)$, the continuous
nest in $c_1$.
\medbreak
For $\N$ countable we have the following natural class:
Let $\alpha$ be a countable ordinal number, index the canonical basis  
of
$\ell_2$ by $\{e_\beta\colon \beta\leq\alpha\}$ and let 
${\bf T}^1(\alpha)$ be the nest algebra in $c_1$ associated to the  
nest
of subspaces $\{N_\beta\colon \beta\leq\alpha\}$ where
$N_\beta=[e_\gamma:\gamma\leq\beta]$.

\medbreak
{THEOREM 6}. No two of the following nest algebras
are isomorphic to each other: $\T^1(\omega)$, $\T^1(2\omega)$, and
$\T^1(\omega^2)$.
\medbreak
If $\alpha\geq\omega^2$ then the intervals of isomorphism
of the $\T^1(\alpha)$'s are at most like those for the $C(\alpha)$'s.
\medbreak
{PROPOSITION 7.} If $\omega^2\leq\alpha\leq\beta<
\alpha^\omega$ then $\T^1(\alpha)\approx\T^1(\beta)$.
\medbreak
The proof of Theorem 5 will consist of two parts, (Lemmas 8 and 9).
The first one shows that $\T^1(\R)$ embeds complementably into
$(\hbox{Alg}\N)^1$, and the second shows that $(\hbox{Alg}\N)^1$
embeds complementably into $\T^1(\R)$. Since $\T^1(\R)\approx
(\sum\oplus\T^1(\R))_1$, the proof follows from the Pe\l czy\'nski
decomposition method \(Pe\).
\medbreak
{LEMMA 8.} If $\N$ is uncountable then $\T^1(\R)$
embeds complementably into $(\hbox{Alg}\N)^1$.
\medbreak
PROOF. It is a consequence of the Similarity Theory \(D\)
that $\N$ is similar to a nest $\tilde{\N}$ with a continuous part.
By \(E\), this one comes from a continuous measure $\mu$ supported
on $\(0,1\)$. Let $P$ be the orthogonal projection onto
$L_2(\(0,1\),\mu)$; then $\Phi(T)=PTP$ sends Alg$\N$ onto a  
continuous
nest. Moreover, $\Phi$ is a projection and also sends 
$(\hbox{Alg}\N)^1$ to a continuous nest in $c_1$.\endpf
\medbreak
LEMMA 9. $(\hbox{Alg}\N)^1$ embeds complementably into
$\T^1(\R)$.
\medbreak
PROOF. The proof uses Erdos' representation Theorem \(E\). For clarity
we will prove it for multiplicity free nests but the proof extends
easily to the general case.

Assume that $\N$ is the standard nest on $L_2(\(0,1\),\mu)$ where
$\mu=\mu_c+\mu_d$ and $\mu_c$ is continuous and $\mu_d$ is discrete
with atoms at $\{d_n\}_n\subset\(0,1\).$

We ``split'' every atom $d_n$ into $d_n^-$ and $d_n^+$ and insert a  
copy
of $\(0,1\)$ in between. $0$ corresponding to $d_n^-$ and 1 to  
$d_n^+$.
More formally, if $I_n=\(0,1\)\times\{n\}$ for every $n$ then
$$\Omega=\left( \(0,1\)\setminus\bigcup_n\{d_n\}\right)\bigcup_n  
I_n$$
with the natural order. It is easy to see that $\Omega$ is then a  
compact
connected space. Define a measure $\nu$ on $\Omega$ by
$\nu=\mu_c$ on $\(0,1\)\setminus\bigcup_n\{d_n\}$ and Lebesgue on
$\bigcup_n I_n$. Then $\nu$ is continuous and we define the standard
continuous nest on $L_2(\Omega,\nu)$, which we denote
Alg$\tilde{\N}$.

To take care of the atoms consider
$x_n=\sqrt{2}\chi_{\(1/2,1\)}$
and $y_n=\sqrt{2}\chi_{\(0,1/2\)}$ supported on $I_n$. Then for every
$n$, $x_n\otimes y_n\in\hbox{Alg}\tilde{\N}$.
($x\otimes y$ denotes the rank-1 map that sends $h\to (h,x)y$).

Let $P$ be the orthogonal projection on $L_2(\Omega,\nu)$ onto
$L_2(\(0,1\)\setminus\bigcup_n\{d_n\},\mu_c)$; $P_x$ the orthogonal
projection onto $\(x_n\)$, and $P_y$ onto $\(y_n\).$ It is clear that
they are orthogonal from each other.

Define $\Phi$ on $B(L_2(\Omega,\nu))$ by
$\Phi(T)=(P+P_x)T(P+P_y).$
Notice that $\Phi$ is a projection and its range is isomorphic to
Alg$\N$. They have the same continuous part: $PTP$; the same atomic
part: $P_xTP_y$; and they interact in the same way. Moreover, $\Phi$
is also defined on $c_1(L_2(\Omega,\nu))$, giving the result.

If $\N$ is not multiplicity free, then represent it as in $\(E\)$,
make the ``enlargement'' on every interval and proceed as before.
\endpf
\bigbreak
The proof of Theorem 6 is more involved. We start with a concrete
representation of $\T^1(\omega)$, (which will be denoted from now on  
by $\T^1$), $\T^1(2\omega)$ and $\T^1(\omega^2)$.
Notice that the last one is isomorphic to
$\T^1\otimes c_1=\overline{\hbox{span}}\{e_{ij}\otimes e_{kl}:i\leq  
j; k,l=1,2,\cdots\}$.

$$\T^1(\omega)=
               \pmatrix{*&*&\cdots\cr
                         &*&\cdots\cr
                         & &\ddots\cr}
,\quad\quad
\T^1(2\omega)=\pmatrix{
               \matrix{*&*&\cdots\cr
                        &*&\cdots\cr
                        & &\ddots\cr}&
               \matrix{*&*&\cdots\cr
                       *&*&\cdots\cr
               \vdots&\vdots&\ddots\cr}\cr
               &
               \matrix{*&*&\cdots\cr
                        &*&\cdots\cr
                        & &\ddots\cr}\cr}
         =\pmatrix{\T^1&c_1\cr
                       &\T^1\cr},\hbox{ and }$$

$$\T^1(\omega^2)=
\pmatrix{\T^1&c_1&c_1&\cdots\cr
            &\T^1&c_1&\cdots\cr
            &   &\T^1&\cdots\cr
            &   &    &\ddots\cr}
\approx
\pmatrix{ c_1&c_1&\cdots\cr
             &c_1&\cdots\cr
             &   &\ddots\cr}
=\T^1\otimes c_1.$$
\medbreak
It is clear that
$\T^1$ embeds complementably into $\T^1(2\omega)$ and this
one embeds complementably into $\T^1(\omega^2)$; however, the
reverse complemented embeddings do not hold (see Lemmas 9 and 10 below).
The key point is the decomposition
$$\T^1(2\omega)\approx\T^1\oplus c_1.$$
\medbreak
LEMMA 10.  $\T^1(2\omega)$ does not embed into $\T^1$.
\medbreak
PROOF. If $\T^1(2\omega)$ embedded into $\T^1$ then
we would have that $c_1$ embeds into $\T^1$. But this is
impossible as we stated in the preliminaries (see \(KP\)).\endpf
\medbreak
The next proposition says that $\T^1(\omega^2)$ does not embed
complementably into $\T^1(2\omega)$.
Nevertheless, since $c_1\subset\T^1(2\omega)$, it does embed.
\medbreak
PROPOSITION 11.
$\T^1(\omega^2)$ does not embed complementably into $\T^1(2\omega)$.
\medbreak
The idea of the proof is that we cannot take away {\it one} $c_1$
from $\T^1\otimes c_1$ in such a way that what we have left is
just {\it one} $\T^1.$

To formalize this we first prove that we can replace any bounded linear operator
$\Phi$ on $\T^1\otimes c_1$ by a multiplier; then we will show
that for this simple type of operator it is not possible to have
$\Phi(\T^1\otimes c_1)\approx c_1$ and
$(I-\Phi)(\T^1\otimes c_1)\approx\T^1.$

A multiplier on $\T^1\otimes c_1$ is a bounded linear operator,
$\Phi$, that satisfies: for every $i\leq j,$
$$\Phi e_{ij}\otimes e_{kl}=\lambda_{ijkl}e_{ij}\otimes e_{kl},$$
for some $\lambda_{ijkl}\in{\bf C}.$

To make the replacement we reduce the problem
from $c_1\otimes c_1$ to $c_1$ and then adapt the steps
of the proof of Theorem 1.

Let $\phi:{\bf N\times N}\to{\bf N}$ be a one-to-one and onto map,
and define $S:c_1\otimes c_1\to c_1$ by
$S(e_{ij}\otimes e_{kl})=e_{\phi(i,k)\phi(j,l)}$.
It is easy to see that $S$ is an isometry onto.
This allows us to define on
$c_1\otimes c_1$ the equivalent maps for $J_{\sigma,\psi}$ and
$K_{\sigma,\psi}$ as follows:
$$\eqalign{\hat{\sigma}(i,j)&
                  =(\phi^{-1}\sigma\phi)(i,j)\quad\hbox{ and}\cr
     \hat{\psi}(k,l)&=(\phi^{-1}\psi\phi)(k,l).\cr}$$
Then $J_{\hat{\sigma},\hat{\psi}}:c_1\otimes c_1\to c_1\otimes c_1$,
and $K_{\hat{\sigma},\hat{\psi}}:c_1\otimes c_1\to c_1\otimes c_1$
are well defined and have the same properties. Moreover, we have
$$S^{-1}J_{\sigma,\psi}=J_{\hat{\sigma},\hat{\psi}}S^{-1}\quad
\hbox{ and }\quad K_{\sigma,\psi}S=S K_{\hat{\sigma},\hat{\psi}}.$$
\medbreak
Let $\S=S({\bf T^1}\otimes c_1).$
Notice that $\S$ is a $*$-diagram; i.e., for $i,j$
fixed,  either $e_{ij}\in\S$ or for every $A\in\S$, $(Ae_j,e_i)=0$.

Let $\pi_1$ be the projection
onto the first coordinate of ${\bf N}\times{\bf N}$ and define
$$r(i)=\pi_1\phi^{-1}(i).$$
It is clear that $e_{ij}\in\S$ if and only if $r(i)\leq r(j).$
Moreover, it is very important to notice that for $i$ fixed there
are infinitely many $j$'s satisfying $r(i)=r(j)$.

\medbreak
LEMMA 12. With the above notation, a necessary and
sufficient condition for
$J_{\sigma,\psi}\S\subset\S$ and $K_{\sigma,\psi}\S\subset\S$ is that
$r(i)\leq r(j)$ if and only if $r(\sigma(i))\leq r(\psi(j)).$
In particular, this is true if $r(i)=r(\sigma(i))=r(\psi(i))$ for
every $i$.
\medbreak
The proof of Lemma 12 follows immediately from the definitions.

\medbreak
LEMMA 13. For every $\epsilon>0$ and
$\tilde{\Phi}:{\bf T^1}\otimes c_1\to {\bf T^1}\otimes c_1$
a bounded operator, there exist $\sigma$ and $\psi$, as in Lemma 12,
and a multiplier $\Phi_1$ on  ${\bf T^1}\otimes c_1$ satisfying
$\|K_{\hat{\sigma},\hat{\psi}}\tilde{\Phi}  
J_{\hat{\sigma},\hat{\psi}}
-\Phi_1\|<\epsilon.$
\medbreak

PROOF OF LEMMA 13. The proof mimics the one of Theorem 1, and we  
solve
it on $\S$ instead. Let $M_n(\S)$, $E_n(\S)$, $F_n(\S)$ and $H_n(\S)$
be as in Theorem 1 with the natural modifications; e.g.,
$M_n(\S)=\hbox{span}\{e_{ij}\in\S:\max\{i,j\}\leq n\}$, etc..

Let $\Phi:\S\to\S$. To simplify notation
whenever we say that $\Phi\approx
\Psi$ we mean that they are arbitrarily close.
\medbreak
STEP 1. Find $\sigma_1$ with $r(i)=r(\sigma_1(i))$ for every $i$,
such that $\Phi_1\approx K_{\sigma_1,\sigma_1}\Phi  
J_{\sigma_1,\sigma_1}$, and it 
satisfies $\Phi_1(M_n(\S))\subset M_n(\S)$
and $\Phi_1(E_n(\S))\subset E_n(\S).$
\medbreak
This is easy. Repeat the proof in Step 1, Theorem 1 but choose
$\sigma(i)$ satisfying
$r(i)=r(\sigma(i))$. This is always possible because
there are infinitely many $j$'s with $r(j)=r(i).$
\medbreak
STEP 2. Find $\sigma_2,\psi_2$, with  
$r(i)=r(\sigma_2(i))=r(\psi_2(i))$
for all $i$, such that $\Phi_2\approx
K_{\sigma_2,\psi_2}\Phi_1 J_{\sigma_2,\psi_2}$ and it satisfies
$\Phi_2(F_n(\S))\subset F_n(\S).$
\medbreak
The proof is very similar to the one of Step 2, Theorem 1. We only  
have to take
$\psi_2(n+1)\in A$ with $r(\sigma_2(n+1))=r(n+1)$ and make sure that
$\{i\colon\psi_2(n)<i\leq N, r(i)=r(n+1)\}$ is ``large enough''
to extract $\sigma_2(n+1)$ from it.

\medbreak
STEP 3. Find $\sigma_3,\psi_3$ as in Lemma 12 such that  
$\Phi_3\approx
K_{\sigma_3,\psi_3}\Phi_2 J_{\sigma_3,\psi_3}$ and it satisfies
$\Phi_3(H_n(\S))\subset H_n(\S).$
\medbreak
The proof is more delicate now. If we repeat the proof of Step 3,
Theorem 1 we may end up with all the elements in $A(n)$ with constant
$r(i)$.

We need to find $\sigma_3,\psi_3$ as in Lemma 12. This means that if
we have chosen
$\sigma_3(i),\psi_3(i)$ for $i\leq n$, then $\sigma_3(n+1)$
and $\psi_3(n+1)$ have the following constrains: If $k\leq n$ then
$$\eqalign{
r(k)\leq r(n+1)&\Longrightarrow r(\sigma_3(k))\leq r(\psi_3(n+1)),
                                              \hbox{ and}\cr
r(k)> r(n+1)&\Longrightarrow r(\sigma_3(k))> r(\psi_3(n+1));\cr}$$
we also have similar conditions for $\sigma_3(n+1)$ in addition to
$r(\sigma_3(n+1))\leq r(\psi_3(n+1)).$ However, we will see that  
there
is a lot of room and we will not worry much.

Let $R_k=\{i:r(i)=k\}$.
We will choose $\sigma_3$ and $\psi_3$ in such a way that if  
$r(i)=r(j)$
then $r(\sigma_3(i))=r(\sigma_3(j))$ and $r(\psi_3(i))=r(\psi_3(j))$.
Therefore, once we choose an element from $R_k$ we must be able to
continue selecting elements from the same $R_k$.

Assume that we have chosen $\sigma_3(i),\psi_3(i)$ for $i\leq n$ and  
that
we have $A(n),B(n)$ subsets of ${\bf N}$ satisfying:
$$\eqalign{
\hbox{card}&(A(n)\bigcap R_{r(\psi_3(i))})=\aleph_0
                                \quad\hbox{ for }i\leq n,\cr
\hbox{card}&\{i:\hbox{ card}(A(n)\bigcap R_i)=\aleph_0\}=\aleph_0,\cr
\hbox{card}&(B(n)\bigcap R_{r(\sigma_3(i))})=\aleph_0
                                \quad\hbox{ for }i\leq n,\cr
\hbox{card}&\{i:\hbox{ card}(B(n)\bigcap  
R_i)=\aleph_0\}=\aleph_0.\cr}
\leqno{(2)}$$

The first and third conditions are the ones that allow us to choose
elements from the previously chosen $R_i$'s; and the others are  
similar
to those of Step 3, Theorem 1.

Suppose we have to choose $\sigma_3(n+1)$ from $R_m$, here $m$ is  
one of the elements of the fourth line of $(2)$. Then let
$\rho\subset B(n)\bigcap R_m$ be such that card($\rho)=N=N(n)$,  where
$N$ is a very large number.

We want to take $A(n+1)\subset A(n)$ satisfying $(2)$, but first we  
find
$A_1(n+1)\subset A(n)$ satisfying $(2)$, and $\rho_1\subset\rho$ very
large such that for $i\leq n$
$$k\in\rho_1\hbox{ and }j\in A_1(n)\bigcap R_{r(\psi_3(i))}
\Longrightarrow \:Q_ne_{kj}\:\leq\epsilon_{n+1}.\leqno{(3)}$$

To check $(3)$ it is enough to do it for only one.  Lemma 2 gives
that for every
$j\in A(n)\bigcap R_{r(\psi_3(1))}$,
card$\{i\in\rho:
\:Q_ne_{kj}\:\geq\epsilon_{n+1}\}$ is very small. Hence, we can find  
a
large $\rho_1\subset\rho$ and $F_1\subset A(n)\bigcap  
R_{r(\psi_3(1))}$
infinite such that if $j\in F_1$ and $k\in\rho_1$ then
$\:Q_ne_{kj}\:\leq\epsilon_{n+1}$. Let
$$A_1(n)=\(A(n)\setminus(A(n)\bigcap R_{r(\sigma_3(1))})\)\bigcup  
F_1.$$
It is clear that this does it. Moreover, it is also clear that we  
have
complete control over
finitely many $R_k$'s. This is the ``room'' we mentioned
before.

So, assume that $A_1(n)$
satisfies $(2)$ and $(3)$. For every $j\in A_1(n)$
there exists $k(j)\in\rho_1$
such that $\:Q_ne_{k(j) j}\:\leq\epsilon_{n+1}.$
Since $\rho_1$ is finite we choose $\sigma_3(n+1)\in\rho$ such that
$A(n+1)=\{j\in A_1(n):j>\max\rho\hbox{ and }\:Q_ne_{\sigma_3(n+1)  
j}\:
\leq\epsilon_{n+1}\}$ satisfies $(2)$. Find $\psi_3(n+1)$ similarly.
\medbreak
STEP 4. Find $\sigma_4,\psi_4$ as in of Lemma 12 such
that $\Phi_4=K_{\sigma_4,\psi_4}\Phi_3 J_{\sigma_4,\psi_4}$ satisfies
$\Phi_4e_{ij}=\lambda_{ij}e_{ij}.$
\medbreak
This is just like Step 4 of Theorem 1.
\medbreak
We finish now the proof of Lemma 13. If
$\tilde{\Phi}\in B(\T^1\otimes c_1)$ then
$S\tilde{\Phi} S^{-1}:\S\to\S$.
Combining Steps 1 through 4 we find $\sigma,\psi$
satisfying the condition of
Lemma 12 and $\tilde{\Phi}_1$ a multiplier such that
$$\:K_{\sigma,\psi}S\tilde{\Phi}  
S^{-1}J_{\sigma,\psi}-\tilde{\Phi}_1\:
\leq\epsilon.$$

Since $K_{\sigma,\psi}S=S K_{\hat{\sigma},\hat{\psi}}$ and
$S^{-1}J_{\sigma,\psi}=J_{\hat{\sigma},\hat{\psi}}S^{-1}$ we obtain  
the
result.\endpf
\medbreak

PROOF OF PROPOSITION 11.
Suppose that $\T^1\otimes c_1\approx\T^1\oplus c_1$.
Then find $\Phi\in B(\T^1\otimes c_1)$ such that $\Phi(\T^1\otimes  
c_1)
\approx c_1$ and $(I-\Phi)(\T^1\otimes c_1)\approx \T^1$. Therefore,
$$\eqalign{
I_{\T^1} &\hbox{ does not factor through }\Phi,\hbox{ and}\cr
I_{c_1}  &\hbox{ does not factor through }I-\Phi.\cr}$$

We will show that this leads to a contradiction.

Find, as in Lemma 13, $\sigma_1,\psi_1$
and a multiplier $\Phi_1$ such that
$\:\Phi_1-K_{\hat{\sigma}_1,\hat{\psi}_1}\Phi
J_{\hat{\sigma}_1,\hat{\psi}_1}\:\leq\epsilon$ and for $i\leq j,$
$\Phi_1e_{ij}\otimes e_{kl}=\lambda_{ijkl}e_{ij}\otimes e_{kl}.$

By a standard diagonal argument we can assume that
$$\eqalign{
\lim_{l\to\infty}\lambda_{ijkl}&=\lambda_{ijk},\cr
\lim_{k\to\infty}\lambda_{ijkl}&=\overline{\lambda}_{ijl},\cr
\lim_{k\to\infty}\lambda_{ijk}&=\lambda_{ij},\cr
\lim_{l\to\infty}\overline{\lambda}_{ijl}&=\overline{\lambda}_{ij}.\cr}$$

CLAIM 1.
$\lambda_{ij}=\overline{\lambda}_{ij}$.
\medbreak
The proof of this is essentially Step 5 of Theorem 1. If for some
$i,j$ we had $\lambda_{ij}\not=\overline{\lambda}_{ij}$, then
looking at the $c_1$ at the $i,j$ position we would find a block
projection with upper triangular part $\lambda_{ij}$ and
lower $\overline{\lambda}_{ij}$. And this would be unbounded.
\medbreak
We can assume moreover that all of the $\lambda_{ij}$
essentially agree; i.e., for some $\lambda\in{\bf C}$,
$\vert \lambda_{ij}-\lambda\vert <
\epsilon({1\over 2^i}+{1\over 2^j}).$
\medbreak

CLAIM 2.
$\vert 1-\lambda\vert\leq 2\epsilon.$
\medbreak
If not, define $J_2:c_1\to\T^1\otimes c_1$ and
$K_2:\T^1\otimes c_1\to c_1$ by
$J_2e_{ij}=e_{ij}\otimes e_{11}$ and $K_2(e_{ij}\otimes  
e_{kl})=e_{ij}$
only if $k=l=1$. Then
$\Phi_2=K_2\Phi_1 J_2:c_1 \to c_1$; also notice that
$I_{c_1}-\Phi_2=K_2(I-\Phi_1)J_2$.

As in Step 5 of Theorem 1 find $J_3,K_3$ in $c_1$ such that
$\Phi_3=K_3\Phi_2 J_3$ satisfies
$\:\Phi_3-\lambda_{11}I_{c_1}\:<\epsilon.$
Hence, if $K=K_3K_2K_{\hat{\sigma}_1,\hat{\psi}_1}$ and
$J=J_{\hat{\sigma}_1,\hat{\psi}_1}J_2J_3$ we have that
$\:K\Phi J-\lambda I_{c_1}\:<2\epsilon.$
Therefore,
$$\:K(I-\Phi)J-(1-\lambda)I_{c_1}\:<2\epsilon.$$
And since $\vert 1-\lambda\vert>2\epsilon$ we have that
$K(I-\Phi)J:c_1\to c_1$ is invertible. This implies that
$I_{c_1}$ factors through $I-\Phi$, a contradiction.
\medbreak
It remains to prove that both Claims 1 and 2 contradict our
assumption.

For this it will be enough to find $\sigma,\psi$ such that for $i\leq  
j$,
$\vert \lambda_{ij\sigma(i)\psi(j)}\vert\geq 1/2.$
Because once we have this
we define $J:\T^1\to\T^1\otimes c_1$ and $K:\T^1\otimes c_1\to\T^1$  
by
$Je_{ij}=e_{ij}\otimes e_{\sigma(i)\psi(j)}$ and $K$ such that
$KJ=I_{\T^1}.$  Then $K\Phi_1J:\T^1\to\T^1$ is a multiplier with big
elements; hence, an adaptation of
the proof of Step 5 of Theorem 1 implies that $I_{\T^1}$ factors
through $\Phi$ giving a contradiction.

The existence of $\sigma,\psi$ will follow from the next claim which
we prove only in the first row of $c_1$'s.
\medbreak
CLAIM 3. Let $B$ be an infinite subset of
${\bf N}$, and for every $k=1,2,\cdots$ let
$$A_k=\{j\in B:\hbox{ card}\{l:\vert \lambda_{1jkl}\vert\geq1/2\}=
\aleph_0\}.$$
Then for some $k$, card$A_k=\aleph_0$.
\medbreak
We know that if we look at the $(c_1)_{1j}$ then there are plenty of
$\lambda_{1jkl}$'s
in the ``lower-right'' corner that satisfy $\vert\lambda_{1jkl}\vert
\approx1$. More specifically, given $l$ large enough,
there exists $k_0$ such that
for $k\geq k_0$, $\lambda_{1jkl}\approx\lambda_{1j}\approx\lambda$.

If Claim 3 were false, for every $k$ fixed, we would have $\vert
\lambda_{1jkl}\vert<1/2$ ``eventually''.
Hence, one can extract arbitrarily
large blocks that look essentially like
$$\pmatrix{\lambda  &\mu_{12} &\cdots&\mu_{1N}\cr
           \lambda  &\lambda  &\cdots&\mu_{2N}\cr
           \vdots   &\vdots   &\ddots&\vdots  \cr
           \lambda  &\lambda  &\cdots&\lambda \cr}.$$
Where $\vert\mu_{ij}\vert<1/2$.
Then, Ramsey's Theorem, used as in Proposition 17
gives us a large submatrix with upper triangular part $\mu$,
$\vert\mu\vert<1/2$ and lower $\lambda$ with
$\vert\lambda\vert>1-2\epsilon$.
Since $N$ is arbitrarily large and the latter matrices are not  
uniformly
bounded we conclude that $\Phi$ is not bounded.
A clear contradiction\endpf
\bigbreak
We finally start the proof of Proposition 7.
For this we need the following
lemma.
\medbreak
LEMMA 14. If $\omega^2\leq\alpha$ then
$\T^1(\alpha)\approx\T^1(\alpha^2).$
\medbreak
PROOF OF PROPOSITION 7. It is clear that Lemma 14 gives that for  
every
$n\in{\bf N}$, $\T^1(\alpha)\approx\T^1(\alpha^n)$. Hence, if
$\alpha\leq\beta<\alpha^\omega$ we can find $n$ such that
$\alpha\leq\beta<\alpha^n.$ Therefore, we have that
$\T^1(\alpha)$ embeds complementably into $\T^1(\beta)$, and this
one into $\T^1(\alpha^n)$.
Since the latter is isomorphic to $\T^1(\alpha)$
and we have that
$\T^1(\alpha)\approx(\sum\oplus\T^1(\alpha))_1$, we finish
the proof using Pe\l czy\'nski's decomposition theorem.\endpf
\medbreak
Before the proof of Lemma 14 we prove this simpler case,
\medbreak
LEMMA 15. If $\omega^2\leq\alpha$ then
$\T^1(\alpha\omega)\approx\T^1(\alpha).$
\medbreak
PROOF. Notice that $(0,\alpha\omega)=\bigcup_{0\leq n<\omega}I_n$  
where
$I_n=(\alpha n,\alpha(n+1)\)$. Therefore, $\ell_2(\alpha\omega)=
(\sum_n\oplus H_n)_2$ where
$H_n=\(e_\eta:\alpha n<\eta\leq\alpha(n+1)\)$.
Taking the ``diagonal'' of this decomposition,
(which is clearly complemented),
we obtain
$$\eqalign{\T^1(\alpha\omega)
 &\approx\left(\sum_n\oplus\T^1(\alpha)\right)_1
                               \oplus\(\T^1\otimes c_1\)\cr
 &\approx\T^1(\alpha)\oplus\T^1(\omega^2)\approx\T^1(\alpha).\endpf\cr}$$
\medbreak
PROOF OF LEMMA 14. Notice that
$(0,\alpha^2)=\bigcup_{0\leq\xi<\alpha}I_\xi$
where $I_\xi=(\alpha\xi,\alpha(\xi+1)\).$
Then repeating a similar argument
as in Lemma 15 we have that
$$\eqalign{\T^1(\alpha^2)
    &\approx\left(\sum_{0\leq\xi<\alpha}\T^1(\alpha)\right)_1\oplus
                                          \(\T^1(\alpha)\otimes  
c_1\)\cr
    &\approx\T^1(\alpha)\oplus\(\T^1(\alpha)\otimes c_1\).\cr}$$
And we finish the proof if we prove

CLAIM 1. If $\omega^2\leq\alpha$ then
$\T^1(\alpha)\otimes c_1\approx\T^1(\alpha).$

Notice first that $(0,\omega\alpha)=\bigcup_{0\leq\xi<\alpha}I_\xi$
where  $I_\xi=(\omega\xi,\omega(\xi+1)\)$. Then we have that
$$T^1(\omega\alpha)\approx\T^1\oplus\(\T^1(\alpha\otimes c_1)\)
                        \approx \T^1(\alpha)\otimes c_1.$$
Hence, it is enough to prove,

CLAIM 2. If $\omega^2\leq\alpha$ then
$\T^1(\alpha)\approx\T^1(\omega\alpha).$

For $\omega^2\leq\alpha<\omega^\omega$
it is enough to take $\alpha=\omega^n$. And
for this case Lemma 15 gives the result.
For $\alpha=\omega^\omega$ we have that
$\omega\omega^\omega=\omega^\omega$.
Actually, if $\gamma\geq\omega$ then  
$\omega\omega^\gamma=\omega^\gamma.$
Hence, for $\alpha\geq\omega^\omega$ we have
$\omega\alpha<\alpha\omega$ (Just take
$\alpha=\omega^\gamma k+\delta$ for some $\delta<\omega^\gamma$ and
$\gamma\geq\omega$). This implies that
$\T^1(\omega\alpha)$ embeds complementably
into $\T^1(\alpha\omega)$; and now Lemma 15 finishes the proof.\endpf


\bigbreak
\noindent{5. PRIMARITY OF OPERATOR SPACES.}
\medbreak
In this section we show that if $\N$ is an
infinite nest in a separable Hilbert space
then Alg$\N$ and $B(H)/{\rm Alg}\N$ are primary.
These are non-commutative versions of theorems proved by
Bourgain \(B\) and M\"uller \(M\"u\).
\medbreak
The technique we use was developed by Bourgain [B] to prove
that $H^\infty$ is primary. It allows to obtain the general
theorem from its finite dimensional version.
\medbreak
To prove that Alg$\N$ is primary we use the decomposition, from  
\(A2\),
$${\rm Alg}\N\approx\left(\sum\oplus\T_n\right)_\infty,$$
where $\T_n$ is the set of all $n\times n$ upper triangular matrices
in $M_n$.
Then we modify slightly the combinatorial argument used by Blower [Bl]
to prove that $B(H)$ is primary. It is the arguments in [Bl] which motivated
our work above.
\medbreak
THEOREM 16. Alg$\N$ is primary.
\medbreak

The proof of the theorem will follow from the finite dimensional
case as Bourgain \(B\) indicated. Since the proof is a
slight modification of \(Bl\) we
use the notation employed there. 
\medbreak

If $\sigma=\{\sigma(1),\cdots,\sigma(n)\}$ and
$\psi=\{\psi(1),\cdots,\psi(n)\}$ are finite subsets of
$\{1,\cdots,N\}$,
define $J_{\sigma,\psi}:M_n\to M_n$ and
$K_{\sigma,\psi}:M_n\to M_n$ by
$$\eqalign{J_{\sigma,\psi}(e_{ij})&=e_{\sigma(i)\psi(j)},\hbox{  
and}\cr
           K_{\sigma,\psi}(e_{ij})&=
     \cases{e_{kl},&if $\sigma(k)=i\hbox{ and }\psi(l)=j$;\cr
              0   ,&otherwise.\cr}\cr}$$

It is easy to see that $J_{\sigma,\psi}$
is an isometric embedding and $K_{\sigma,\psi}J_{\sigma,\psi}=I$,
where $I$ is the identity of $M_n.$

Moreover, if
$$\sigma(1)<\psi(1)<\sigma(2)<\psi(2)
<\cdots<\sigma(n)<\psi(n),\leqno{(4)}$$
then $K_{\sigma,\psi}(\T_N)\equiv\T_n$.
\medbreak
Then we obtain a proposition similar to the one in \(Bl\).
\medbreak
PROPOSITION 17. Given $n, \epsilon>0$ and $K<\infty$
there exists $N_0$ such that if $N>N_0$ and $T\in B(\T_N,\T_N)$
with $\|T\|\leq K$, then there exist subsets $\sigma$ and $\psi$
of $\{1,\cdots,N\}$ of cardinality $n$ such that
$\sigma(1)<\psi(2)<\cdots<\sigma(n)<\psi(n)$ and a constant $\lambda$  
such
that
$$\|K_{\sigma,\psi}TJ_{\sigma,\psi}-\lambda I_n\|\leq\epsilon.$$
Thus, one of $K_{\sigma,\psi}TJ_{\sigma,\psi}$ and
$K_{\sigma,\psi}(I_n-T)J_{\sigma,\psi}$ is invertible.
\medbreak
REMARK. The previous proposition was proved by Blower
[Bl] without the assumption that $\sigma$ and $\psi$
satisfy (4).  In fact, he proved it for 
$\max\sigma<\min\psi$. However, since we need to preserve
the triangular structure we modify the argument from \(Bl\)
to obtain the desired result.
\medbreak
SKETCH OF THE PROOF. Just as in \(Bl\) find a large
$\sigma_1\subset\{1,\cdots,N\}$ and $\lambda\in{\bf C}$ such that if
$i<j$ are in $\sigma_1$ then
$\vert (T(e_{ij})_{ij}-c \vert <\epsilon n^{-6}/4.$

The goal now is to find a large $\tilde{\sigma}\subset\sigma_1$
such that if $i<j$ and $k<l$ are in $\tilde{\sigma}$ and
satisfy $(i,j)\not=(k,l)$ then $\vert(Te_{ij})_{kl}\vert<\delta$,
where $\delta=\epsilon n^{-6}/4$. Then
$\sigma=\{\tilde{\sigma}({2i-1})\}_{i=1}^n$ and
$\psi=\{\tilde{\sigma}({2i})\}_{i=1}^n$ will do it.

To find $\tilde{\sigma}$ we find first, as in \(Bl\), a large
$\sigma_2\subset\sigma_1$ such that if $i<k<j<l$ are in $\sigma_2$
then $\vert(Te_{ij})_{kl}\vert<\delta$. Then find a large
$\sigma_3\subset\sigma_2$ such that if $k<i<j<l$ are in $\sigma_3$
then $\vert(Te_{ij})_{kl}\vert<\delta.$ After this find a large
$\sigma_4\subset\sigma_3$ such that if $i=k<j<l$ are in $\sigma_4$
one has $\vert(Te_{ij})_{kl}\vert<\delta.$  Proceeding in this way
we finish. \endpf

\medbreak
LEMMA 18. Given $n\in {\bf N}$, $\epsilon>0$,
there exists an $N'(n,\epsilon)$ such that if $N>N'(n,\epsilon)$
and $E$ is an $n$-dimensional subspace of $\T_N$ then there
exists a subspace $F$ of $\T_N$ and a block projection
$q$, satisfying (2), from $\T_N$ onto $F$ such that
$\|qx\|\leq\epsilon\|x\|$ for $x\in E$.
\medbreak
PROOF. It is enough to show that if $x\in\T_n$,
$\|x\|=1$ then we can find $q$, a large
block projection that respects triangularity, such that
$\|q(x)\|\leq\epsilon$.  Then take an $\epsilon$-net
of the unit sphere of $E$, $\{x_i\}_{i=1}^M$. Find $q_1$ a large
block projection satisfying (4) such that $\|q_1(x_1)\|<\epsilon$;
after this find $q_2$ a block projection contained in the range of  
$q_1$
such that $\|q_2(q_1x_2)\|<\epsilon$. Proceeding in this way we get
$q=q_M\cdots q_2q_1$; and $q$ does it.

To check the previous remark let $x\in\T_N$, $\|x\|=1$,
$\delta>0$ (to be fixed latter), and let
$$\{i,j\}\ \hbox{ is {\it bad} if }\ i<j\ \hbox{ and }\
\vert x_{ij}\vert \geq\delta.$$
Ramsey's Theorem gives us a large monochromatic subset $\rho$. If  
$\rho$
were {\it bad} we would set $i=\min\rho$ and then
$$\|x(e_i)\|^2\geq\sum_{j\in\rho\setminus\{i\}}\vert x_{ji}\vert^2
              \geq(\vert\rho\vert-1)\delta^2.$$
Since, $\|x(e_i)\|\leq1$, the right choice of $\delta$ would give us
a contradiction; hence, $\rho$ is {\it good}.\endpf
\medbreak
We conclude with some comments on the proof
of the following proposition,
\medbreak
PROPOSITION 18. $B(H)/{\rm Alg}\N$ is primary.
\medbreak
The proof is similar to the one for Alg$\N$.
We have an isomorphic representation for $B(H)/{Alg}\N$
similar to the one for Alg$\N$; i.e.,
$$B(H)/\hbox{Alg}\N\approx
\left(\sum_{n=1}^\infty\oplus M_n/\T_n\right)_\infty,$$
where $\T_n$ is the algebra of all upper triangular $n\times n$  
matrices.
This follows from the proof of the main result in \(A2\).
The reason for this is that the maps
$\phi_n:\hbox{Alg}\N\to\A_n$  and $\psi_n:\A_n\to\hbox{Alg}\N$
from \(PPW\) are restrictions from
$\tilde{\phi}_n:B(H)\to\B_n$ and $\tilde{\psi}_n:\B_n\to B(H)$,
where $\B_n$ is the enveloping algebra
of $\A_n$. Therefore, this induces maps
$\phi'_n:B(H)/\hbox{Alg}\N\to\B_n/\A_n$ and
$\psi'_n:\B_n/\A_n\to B(H)/\hbox{Alg}\N$
with the right properties.

Once we have the decomposition, the combinatorial argument
is essentially the same.
\bigbreak
\noindent{6. APPLICATIONS AND OPEN QUESTIONS.}
\medbreak
We conclude this paper with some applications and open questions. The
first one is that the nest algebras in the trace class have bases.

It is well known that $c_1$ has a basis. We take the elements along  
the
``shell'' decomposition $\{F_n\}_n$ (See the proof of Theorem 1 for  
the
notation) of $c_1$; i.e,
$$e_{11},\ \ e_{12},e_{22},e_{21},\ \
e_{13},e_{23},e_{33},e_{32},e_{31},\cdots\ \hbox{ etc.}.$$
Therefore, if $\S\subset c_1$ is a $*$-diagram,
(i.e., for $i,j$ fixed either $e_{ij}\in \S$
or for every $A\in\S$ we have $(Ae_j,e_i)=0$), it has a basis.
This is so because we are just
taking a subsequence of the basis.

This is the principle we use to prove:
\medbreak
THEOREM 19. The nest algebras in $c_1$ have bases.
\medbreak
PROOF. We divide the proof in two cases:
If the nest is uncountable
we use Theorem 7 and the following representation.
Index the basis of $\ell_2$ by the rational numbers in $\(0,1\)$,
$\{e_r\}_r$. Then for $0\leq t\leq1$, let
$N_t=\(e_r:r\leq t\)$ and $N_t^-=\(e_r:r<t\)$.  It is clear that
$\N=\{N_t,N_t^-\}_{0\leq t\leq1}$
is an uncountable nest algebra. Therefore,
if $\phi:{\bf N}\to {\bf Q}\bigcap\(0,1\)$
is a one-to-one and onto map and
$U:\ell_2\to\ell_2({\bf Q}\bigcap\(0,1\))$
is defined by $Ue_i=e_{\phi(i)}$,
we have that $U^{-1}(\hbox{Alg}\N)^1 U$
is a $*$-diagram in $c_1$. Therefore, it has a basis.

If $\N$ is countable we use the representation theorem from $\(E\)$.
For simplicity we do it only for the multiplicity free case, although
the proof for the general case is basically the same.

Assume then that our nest is the standard nest on $L_2(\(0,1\),\mu)$  
for
some measure $\mu$. Since $\N$ is countable we have that $\mu$ is
totally atomic. Therefore, if we do a similar construction as above,  
we
see that it is unitarily equivalent to a $*$-diagram and
therefore it has a basis. \endpf
\medbreak
REMARK. Most of the results of this paper work in the space of
compact operators with some minor notational changes.
In particular, Theorem 5, part of Theorem 6, Proposition 7 and
Theorem 19 all hold in $K(H)$.
\medbreak
We conclude this three questions.
\medbreak
QUESTION 1. If $\N$ is a countable nest, does $(\hbox{Alg}\N)^1$
correspond to some $\T^1(\alpha)$?
\medbreak
This question is motivated by the analogy between the classification
of the space of continuous functions on compact metric spaces and the
nest algebras in $c_1$. Recall that $(\N,\leq)$ is a compact space that
can be taken inside $\(0,1\)$. If $\N$ is uncountable, then
$C(\N)\approx C(\(0,1\))$ and $(\hbox{Alg}\N)^1\approx\T^1(\R)$.
If $\N$ is countable, then $C(\N)\approx C(\omega^\alpha)$ where
$\alpha$ is the smallest ordinal number for which
$\N^{(\alpha)}$, (the $\alpha$th derived set of $\N$), is finite.
We cannot reproduce the previous result exactly, because for the
nest algebra case it matters if the limit points are one-sided or
two-sided. For example, if we take $A_1=\{1/2-1/n\}_n$ and
$A_2=\{1/2\pm 1/n\}_n$ then $C(A_1)\approx C(A_2)$ but
$\T^1(A_1)\approx\T^1(\omega)$ and $\T^1(A_2)\approx\T^1(2\omega).$
Nevertheless, they correspond to some $\T^1(\alpha)$. The problem
seems to be at the limit points; i.e., if $\N^{(\alpha)}$ is
finite and $\alpha$ is a limit point, say $\omega$.
\medbreak
QUESTION 2. Are there uncountable many non-isomorphic $\T^1(\alpha)$'s?
\medbreak
In particular we are asking if  
$\T^1(\alpha)\approx\T^1(\alpha^\omega)$.
A first step to question 2 is
\medbreak
QUESTION 3. is $\T^1(\omega^2)\approx\T^1(\omega^\omega)$?

\bigbreak
\centerline{REFERENCES}
\medbreak
\item{\(AL\)} J.~Arazy and J.~Lindenstrauss, {\sl Some linear  
topological
properties of the space $c_p$ on operators on Hilbert spaces},
Compositio Math. {\bf 30} (1975), 81--111.

\item{\(Ar1\)} J.~Arazy, {\sl On subspaces of $c_p$ which contain  
$c_p$},
Composition Math. {\bf 41} (1980), 297--336.

\item{\(Ar2\)} J.~Arazy, {\sl Basic sequences, embeddings, and  
uniqueness
of the symmetric structure in unitary spaces}, J. Funct. Anal.
{\bf 40} (3) (1981), 302--340.

\item{\(A1\)} A.~Arias,
{\sl Some properties of non-commutative $H^1$-spaces},
Proc. Amer. Math. Soc. {\bf 112} (2) (1991), 465--472.

\item{\(A2\)} A.~Arias, {\sl All the nest algebras are isomorphic},
Proc. Amer. Math. Soc. {\bf 115} (1) (1992), 85--88.

\item{\(BP\)} C.~Bessaga and A.~Pe\l czy\'nski, {\sl Spaces of  
continuous
functions (IV)}, Studia Math. {\bf 19} (1960), 53--62.

\item{\(Bl\)} G.~Blower, {\sl The space $B(H)$ is primary},
Bull. London Math Soc. {\bf 22} (1990), 176--182.

\item{\(Bo\)} B.~Bollob\'as, {\sl Combinatorics}, (Cambridge
University Press, 1986).

\item{\(B\)} J.~Bourgain, {\sl On the primarity of  
$H^\infty$-spaces},
Israel J. Math {\bf 45} (1983), 329--336.

\item{\(D\)} K.~Davidson,
{\sl Similarity and compact perturbations of nest
algebras}, J.~Reine Angew. Math. {\bf 348} (1984), 286--294.

\item{\(E\)} J.~Erdos, {\sl Unitary invariant for nests}, Pacific  
J.~Math
{\bf 23} (1967), 229--256.

\item{\(FAM\)} T.~Fall, W.~Arveson and P.~Muhly, {\sl Perturbations  
of
nest algebras}, J.~Operator Theory {\bf 1} (1979), 137--150.

\item{\(GK\)} I.~Gohberg and M.~Krein, Theory and applications of
Volterra operators in Hilbert space, {\sl ``Nauka'', Moskow, 1967};
English trans., Trans. Math. Monographs, {\bf 24}, AMS,
Providence, RI, 1970.

\item{\(KS\)} R.~Kadison and I.~Singer, {\sl Triangular operator
algebras}, Amer. J. Math {\bf 82} (1960), 227--259.

\item{\(KP\)} S.~Kwapien and A.~Pe\l czy\'nski, {\sl The main  
triangular
projection in matrix spaces and its applications}, Studia Math. {\bf 34} (1970),
43--68.

\item{\(L\)} D.~Larson,
{\sl Nest algebras and similarity transformations},
Ann. Math. {\bf 121} (1988), 409--427.

\item{\(M\)} V.~I.~Macaev,
{\sl Volterra operators produced by perturbations
of selfadjoint operators}, Soviet Math. Dokl. {\bf 2} (1961),  
1013--1016.

\item{\(M\"u\)} P.~M\"uller, {\sl On projections in $H^1$ and BMO},
Studia Math. (1988), 145--158.

\item{\(PPW\)} V.~Paulsen, S.~Power and J.~Ward, {\sl  
Semidiscreteness
and dilations of nest algebras}, J. Funct. Anal. {\bf 80} (1) (1988),
76--87.

\item{\(Pe\)} A.~Pe\l czy\'nski,
{\sl Projections in certain Banach spaces},
Studia Math. {\bf 19} (1960), 209--228.

\item{\(P\)} S.~Power, {\sl Analysis in nest algebras}, Surveys of  
Recent
Results in Operator Theory, J.~Conway (ed.), Pitman Research Notes in
Mathematics, Longman.

\item{\(R\)} J.~R.~Ringrose, {\sl On some algebras of operators},  
Proc.
London Math. Soc. {\bf 15} (3) (1965), 61--83.
\bigbreak

\noindent Current Address:

\noindent Division of Mathematics, Computer Science and Statistics,

\noindent The University of Texas at San Antonio,

\noindent San Antonio, Texas 78249.

\bye